\documentstyle[11pt]{article}
\newcommand{\Cos}{\mathop{\rm Cos}\nolimits}
\newcommand{\Sin}{\mathop{\rm Sin}\nolimits}
\newcommand{\dn}{\mathop{\rm dn}\nolimits}
\newcommand{\R}{\mathop{\rm Re}\nolimits}

\newcommand{\iiqq}{J_\nu^{(2)}((1-q^2)q^{-1/2}zs;q^2)}
\newcommand{\iqo}{\ddag J_\nu^{(1)}((1-q^2)zs;q^2)\ddag}
\newcommand{\iio}{J_0^{(2)}((1-q^2)q^{-1/2}zs;q^2)}
\newcommand{\io}{\ddag J_1^{(1)}((1-q^2)zs;q^2)\ddag}
\newcommand{\ioq}{\ddag J_1^{(1)}((1-q^2)q^{-1}zs;q^2)\ddag}
\newcommand{\iioq}{J_0^{(2)}((1-q^2)q^{-3/2}zs;q^2)}
\newcommand{\iioqq}{J_0^{(2)}((1-q^2)q^{-5/2}zs;q^2)}
\newcommand{\iiqo}{\ddag J_\nu^{(2)}((1-q^2)zs;q^2)\ddag}
\newcommand{\qpop}{\frac{1-q^2}2}
\newcommand{\pros}{\frac{2\partial_s}{1+q}}
\newcommand{\proz}{\frac{2\partial_z}{1+q}}
\newcommand{\bg}{\phantom._0\Phi_1(-;0;q,\frac{1-q^2}2zs)}

\newcommand{\ibg}{\phantom._0\Phi_1(-;0;q,i\frac{1-q^2}2zs)}

\newcommand{\exo}{\ddag e_q(\frac{1-q^2}2(zs))\ddag }
\newcommand{\exoa}{\ddag e_q(\frac{1-q^2}2azs)\ddag }

\newcommand{\exoq}{\ddag e_q(\frac{1-q^2}2q^{-1}zs)\ddag }

\newcommand{\Exo}{\ddag E_q(\frac{1-q^2}2zs)\ddag }
\newcommand{\exoi}{\ddag e_q(i\frac{1-q^2}2zs)\ddag }
\newcommand{\exoiq}{\ddag e_q(i\frac{1-q^2}2q^{-1}zs)\ddag}
\newcommand{\Eexiqq}{E_q(i\frac{1-q^2}2q^{-2}zs)}
\newcommand{\Eexiq}{E_q(i\frac{1-q^2}2q^{-1}zs)}
\newcommand{\Eexi}{E_q(i\frac{1-q^2}2zs)}
\newcommand{\Exoa}{\ddag E_q(\frac{1-q^2}2azs)\ddag }

\newcommand{\iExo}{\ddag E_q(i\frac{1-q^2}2zs)\ddag }

\newcommand{\iExop}{\ddag E_q(i\frac{1-q^2}2qzs)\ddag }
\newcommand{\Eex}{E_q(\frac{1-q^2}2zs)}

\newcommand{\Eexq}{E_q(\frac{1-q^2}2q^{-1}zs)}
\newcommand{\Eexqq}{E_q(\frac{1-q^2}2q^{-2}zs)}
\newcommand{\bbg}{\phantom._0\Phi_3(-;0,0,q^{2\nu+2};q^2,
-(\frac{1-q^2}2q^{2(\nu+1)-1/2}zs)^2)}

\newcommand{\G}{\Gamma}
\newcommand{\p}{\partial}

\newcommand{\mbes}{I_\nu^{(1)}((1-q^2)s;q^2)}

\newcommand{\mbess}{I_\nu^{(2)}((1-q^2)s;q^2)}
\newcommand{\makd}{K_\nu^{(1)}((1-q^2)s;q^2)}
\newcommand{\makdo}{K_\nu^{(2)}((1-q^2)s;q^2)}
\newcommand{\makdj}{K_\nu^{(j)}((1-q^2)s;q^2)}

\newtheorem{predl}{Proposition}[section]
\newtheorem{defi}{Definition}[section]
\newtheorem{rem}{Remark}[section]
\newtheorem{cor}{Corollary}[section]
\newtheorem{lem}{Lemma}[section]

\newcommand{\beq}[1]{\begin{equation}\label{#1}}
\newcommand{\eq}{\end{equation}}

\begin{document}

\begin{flushright}
q-alg/96xxxxx \\
 ITEP-TH-7/96\\
MPI-96-54\\
\end{flushright}
\vspace{10mm}
\begin{center}
{\Large \bf  Jackson Integral Representations of
 Modified $q$-Bessel
Functions and  $q$-Bessel-Macdonald Functions}\\ \vspace{5mm}
M.A.Olshanetsky \footnote{Supported in part by RFFI-96-02-18046,
 INTAS-94-2317 and NWO grants} \\ ITEP,117259,
Moscow, Russia\\e-mail olshanez@vxdesy.desy.de\\
\vspace{5mm}
V.-B.K.Rogov \footnote{Supported by ISF grant JC 7100}\\ MIIT,101475,
 Moscow, Russia\\
e-mail m10106@sucemi.bitnet\\
February 1996\\
\end{center}
\begin{abstract}
The $q$ analog of Modified Bessel functions and Bessel-Macdonald
functions,
were defined in our previous work (q-alg/950913) as general solutions
of a second order difference equations.
Here we present a collection of their representations
by the Jackson q-integral.
\end{abstract}

\section{Introduction}
\setcounter{equation}{0}
In this work we continue our investigations of
  Modified $q$-Bessel
Functions (qMBF) and  $q$-Bessel-Macdonald Functions (qBMF) introduced
in \cite{OR1}. Though these functions can be easily expressed through the
$q$-Bessel functions (qBF) of Jackson
 \cite{Ja} , as their classical counterparts,
 they are important in their own turn. The classical Bessel functions
 and the qBF arises as matrix elements
in irreducible representations of the group of motion of Euclidean spaces
\cite{Vi} or its quantum deformations \cite{VK1,FV,Kol}. Similarly,
the classical MBF and BMF are related
 to the Whittaker model of irreducible unitary representations of real
simple split Lie group $G$, when $G=SL(2,{\bf R})$ \cite{J,Sh,Ko,Ha,STS}.
 Modified Bessel Functions generate a basis in common
 eigenspace of Laplace
operators in the coordinate systems, related to the Iwaswa decomposition of
 $SL(2,{\bf R})$, while the  Bessel-Macdonald Functions are bounded
eigenvectors.
They are special matrix elements in the irreducible principle series
 representations.
This fact means that BMF have integral representations coming from the
  invariant
Hermitian measure on sections of  line bundles over the flag variety,
 where the irreducible unitary representations
are realized (Borel-Weyl theorem).

The main object of present work is the analogous
integral representations for qMBF and  some
other integral representations for both kind of functions,
 which have well-known
classical forms. Though there is a clear group theoretical
 interpretation of these representations
our derivation is pure analytic, beyond the group methods.
In fact qBMF, we consider here, plays the role
of the Whittaker function for
$U_q(SL(2,{\bf C}))$ \cite{OR2}. The $q$-integral we use 
  very similar to the Fourier transform on the
quantum line and plane \cite{CZ}. Such ingredients as
the horospheric projection,  Whittaker vectors,
 Harish-Chandra function which occur in the harmonic
analysis on quantum Lobachevsky space can be extracted from
 our integrals. They correspond
to the analogous elements in the spherical model of unitary
representation of
$U_q(SU(1,1))$ \cite{VK}, where the Jackson integral was used
to describe the zonal spherical functions.

Another applications of integral representations of qMBF
 is  Poisson kernels
for quantum Lobachevsky space , or more generally, quantum hyperboloids.
The details of both types of interpretations will be published
 elsewhere.

\section{Some preliminary relations}
\setcounter{equation}{0}
 Jackson $q$-integral is determined as the map from an algebra of the
 functions of one variable
into a set of the number series
$$
\int_{-1}^1f(x)d_qx = (1-q)\sum_{m=0}^\infty q^m[f(q^m)+f(-q^m)],
$$
$$
\int_0^\infty f(x)d_qx=(1-q)\sum_{m=-\infty}^\infty q^mf(q^m),
$$
$$
\int_{-\infty}^\infty f(x)d_qx = (1-q)\sum_{-\infty}^\infty q^m
[f(q^m)+f(-q^m)].
$$
Define the difference operator
\begin{equation}
\label{2.1}
\partial_xf(x)=\frac{x^{-1}}{1-q}[f(x)-f(qx)].
\end{equation}
The following formulas of the $q$-integration by parts are valid
\begin{equation}
\label{2.2}
\int_{-1}^1\phi(x)\partial_x\psi(x)d_qx=
\phi(1)\psi(1)-\phi(-1)\psi(-1)-
\int_{-1}^1\partial_x\phi(x)\psi(qx)d_qx,
\end{equation}
\begin{equation}\label{2.3}
\int_0^\infty\phi(x)\partial_x\psi(x)d_qx=
\lim_{m\to\infty}[\phi(q^{-m})\psi(q^{-m})-
\phi(q^m)\psi(q^m)]-
\int_0^\infty\partial_x\phi(x)\psi(qx)d_qx.
\end{equation}
\begin{equation}\label{2.4}
\int_{-\infty}^\infty\phi(x)\partial_x\psi(x)d_qx=
\lim_{m\to\infty}[\phi(q^{-m})\psi(q^{-m})+
\phi(-q^{-m})\psi(-q^{-m})]-
\int_{-\infty}^\infty\partial_x\phi(x)\psi(qx)d_qx.
\end{equation}
The last two expressions imply the regularization of the imprope
 integrals.

Let $z$ and $s$ be  noncommuting elements and
\begin{equation}
\label{2.5}
zs=qsz.
\end{equation}

Consider the function
\begin{equation}
\label{2.6}
f(x)=\sum_ra_rx^r.
\end{equation}
\bigskip
\noindent
{\it The rule of $q$-integration} in the noncommutative case is
$$
\int f(zs)d_qs=\int\sum_ra_r(zs)^rd_qs=
\int\sum_ra_rq^{-\frac{r(r-1)}2}z^rs^rd_qs,
$$
$$
\int d_qzf(zs)=\int d_qz\sum_ra_r(zs)^r=
\int d_qz\sum_ra_rq^{-\frac{r(r-1)}2}z^rs^r.
$$
 Define the following transformation $\ddag f\ddag$ for functions
 $f$ depending
on the noncommutative variables $s$ and $z$ (\ref{2.5}).
If we have the function which has the form (\ref{2.6}) 
and all monoms are ordered we will write
$$
f(zs)=\sum_ra_r(zs)^r\rightarrow \ddag f(zs)\ddag
=\sum_ra_rz^rs^r.
$$
\begin{defi}\label{d2.1}
The function $f(z)$ is absolutely  $q$-integrable if the series
$$
\sum_{m=-\infty}^\infty q^mf(q^m)
$$
converges absolutely.
\end{defi}
It means, in particular, that
$$
\lim_{m\to\pm\infty}q^m|f(q^m)|=0
$$

Let \cite{GR}
$$
(a,q)_n=\left\{
\begin{array}{lcl}
1&for&n=0\\
(1-a)(1-aq)\ldots(1-aq^{n-1})&for&n\ge1,\\
\end{array}
\right.
$$
\begin{equation}
\label{2.7}
(a,q)_\infty=\lim_{n\to\infty}(a,q)_n,\qquad
(a_1,\ldots,a_k,q)_\infty=(a_1,q)_\infty\ldots(a_k,q)_\infty.
\end{equation}
and $$
\Gamma_q(\nu)=\frac{(q,q)_{\infty}}{(q^{\nu},q)_\infty}(1-q)^{1-\nu}
$$
is the $q$-gamma function.
Consider the $q$-exponentials
\begin{equation}
\label{2.8}
e_q(\frac{1-q^2}2x)=\sum_{n=0}^\infty\frac{(1-q^2)^nx^n}{(q,q)_n2^n},
\qquad|x|<\frac2{1-q^2},
\end{equation}
\begin{equation}
\label{2.9}
E_q(\frac{1-q^2}2x)=\sum_{n=0}^\infty
\frac{q^{\frac{n(n-1)}2}(1-q^2)^nx^n}{(q,q)_n2^n},
\end{equation}
Note that
\beq{2.9a}
e_q(\frac{1-q^2}2x)=\frac1{(\frac{1-q^2}2x,q)_\infty},~~
E_q(\frac{1-q^2}2x)=(-\frac{1-q^2}2x,q)_\infty.
\eq
and
\beq{2.9b}
e_q(q)E_q(q^{-1})=\frac 1{(q,q)_\infty}\sum_{k=0}^\infty
\frac
{q^{\frac{k(k+1)}{2}}}
{(q,q)_k}
\eq
 Consider also the basic hypergeometric series
$$
\phantom._0\Phi_1(-;0;q,\frac{1-q^2}2x)=
\sum_{n=0}^\infty\frac{q^{n(n-1)}(1-q^2)^nx^n}{(q,q)_n2^n}.
$$
It follows immediately
\begin{predl}\label{p2.1}
\begin{equation}
\label{2.10}
\bg=\Exo,
\end{equation}
\begin{equation}
\label{2.11}
\Eex=\exo.
\end{equation}
\end{predl}
\bigskip
Now consider the $q$-Bessel function \cite{Ja}
\begin{equation}
\label{2.12}
J_\nu^{(1)}((1-q^2)x;q^2)=\frac1{\Gamma_{q^2}(\nu+1)}
\sum_{n=0}^\infty
\frac{(-1)^n(1-q^2)^{2n}x^{\nu+2n}}
{(q^2,q^2)_n(q^{2\nu+2},q^2)_n2^{\nu+2n}},
\end{equation}
\begin{equation}
\label{2.13}
J_\nu^{(2)}((1-q^2)x;q^2)=\frac1{\Gamma_{q^2}(\nu+1)}
\sum_{n=0}^\infty
\frac{(-1)^nq^{2n(\nu+n)}(1-q^2)^{2n}x^{\nu+2n}}
{(q^2,q^2)_n(q^{2\nu+2},q^2)_n2^{\nu+2n}},
\end{equation}
and the basic hypergeometric series \cite{GR}
\begin{equation}
\label{2.14}
\phantom._0\Phi_3(-;0,0,q^{2\nu+2};q^2,-(\qpop q^{2\nu+2}x)^2)=
\sum_{n=0}^\infty
\frac{(-1)^nq^{4n(\nu+n)}(1-q^2)^{2n}x^{2n}}
{(q^2,q^2)_n(q^{2\nu+2},q^2)_n2^{2n}}.
\end{equation}

As a generalization of (\ref{2.5}) we assume for any $\nu$ and 
$\mu$
\begin{equation}
\label{2.15}
z^\nu s^\mu=q^{\nu\mu}s^\mu z^\nu.
\end{equation}
\begin{predl}\label{p2.2}
$$
\frac1{\Gamma_{q^2}(\nu+1)}\left(q^{-1/2}\frac{zs}2\right)^\nu\bbg=
$$
\begin{equation}\label{2.16}
=q^{-\frac{\nu^2}2}\iiqo,
\end{equation}
\begin{equation}
\iiqq=q^{-\frac{\nu^2}2}\iqo.
\end{equation}\label{2.17}
\end{predl}

{\bf Proof.} These formulas follow from (\ref{2.15}) and
(\ref{2.12})-(\ref{2.14}).\rule{5pt}{5pt}
\bigskip

In the noncommutative case we define the difference operator similarly
 to (\ref{2.1})
$$
\p_sf(zs)=\frac{s^{-1}}{1-q}[f(zs)-f(qzs)].
$$
Then
$$
\frac1{1-q}[f(zs)-f(qzs)]s^{-1}=\p_sf(q^{-1}zs).
$$
As it follows from (\ref{2.1}) and (\ref{2.8}), (\ref{2.9}) for
arbitrary  $a$
\begin{equation}
\label{2.18}
\proz\exoa=a\exoa s,
\end{equation}
\begin{equation}
\label{2.19}
\pros\exoa=aqz\ddag e_q(\qpop aq(zs))\ddag,
\end{equation}
\begin{equation}
\label{2.20}
\proz\Exoa=a\ddag E_q(\qpop aq(zs))\ddag s,
\end{equation}
\begin{equation}
\label{2.21}
\pros\Exoa=aqz\ddag E_q(\qpop aq^2(zs))\ddag.
\end{equation}
Similarly from (\ref{2.1}) and (\ref{2.12}), (\ref{2.13})
\begin{equation}
\label{2.22}
\proz[(z/2)^{-\nu}\ddag J_\nu^{(1)}((1-q^2)azs;q^2)\ddag]=
-a(z/2)^{-\nu}\ddag J_{\nu+1}^{(1)}((1-q^2)azs;q^2)\ddag s,
\end{equation}
\begin{equation}
\label{2.23}
\pros[\ddag J_\nu^{(1)}((1-q^2)azs;q^2)\ddag(s/2)^{-\nu}]=
-aqz\ddag J_{\nu+1}^{(1)}((1-q^2)aqzs;q^2)\ddag(s/2)^{-\nu},
\end{equation}
\begin{equation}
\label{2.24}
\proz[(z/2)^\nu\ddag J_\nu^{(1)}((1-q^2)azs;q^2)\ddag]=
a(z/2)^\nu\ddag J_{\nu-1}^{(1)}((1-q^2)azs;q^2)\ddag s,
\end{equation}
\begin{equation}
\label{2.25}
\proz[(z/2)^{-\nu}\ddag J_\nu^{(2)}((1-q^2)azs;q^2)\ddag]=
-aq^{\nu+1}(z/2)^{-\nu}\ddag J_{\nu+1}^{(2)}((1-q^2)aqzs;q^2)\ddag s,
\end{equation}
\begin{equation}
\label{2.26}
\pros[\ddag J_\nu^{(2)}((1-q^2)azs;q^2)\ddag(s/2)^{-\nu}]=
-aq^{\nu+2}z\ddag J_{\nu+1}^{(2)}((1-q^2)aq^2zs;q^2)\ddag(s/2)^{-\nu},
\end{equation}
\begin{equation}
\label{2.27}
\proz[(z/2)^\nu\ddag J_\nu^{(2)}((1-q^2)azs;q^2)\ddag]=
aq^{-\nu+1}(z/2)^\nu\ddag J_{\nu-1}^{(2)}((1-q^2)aqzs;q^2)\ddag s,
\end{equation}

Consider the complete elliptic integrals:
$$
{\bf K}(k)=\int_0^{\pi/2}\frac{d\alpha}{\sqrt{1-k^2\sin^2\alpha}},
\qquad {\bf K}'(k)=\int_0^{\pi/2}\frac{d\alpha}
{\sqrt{\cos^2\alpha+k^2\sin^2\alpha}}
$$
with
\begin{equation}\label{2.30}
\ln q=-\frac{\pi{\bf K}'(k)}{{\bf K}(k)}.
\end{equation}
Then \cite{BMP}(5.3.6.1)
\begin{equation}\label{2.31}
{\bf Q}_\nu=(1-q)\sum_{m=-\infty}^\infty\frac1
{q^{m-\nu+1/2}+q^{-m+\nu-1/2}}=\frac{1-q}{\pi}{\bf K}(k)
\dn[\frac{2\ln q^{-\nu+1/2}}{\pi}{\bf K}'(k)],
\end{equation}
where $\dn u=\sqrt{1-k^2\sin^2\phi},\quad u=\int_0^\phi\frac{d\alpha}
{\sqrt{1-k^2\sin^2\alpha}}$, 
($\dn u$ is the Jacobi elliptic function).
If $u=0$ then $\phi=0$, and $\dn u=1$.
\begin{predl}\label{p2.3}
For an arbitrary $\nu$
\begin{equation}\label{2.32}
\lim_{q\to 1-0}{\bf Q}_\nu=\frac\pi2
\end{equation}
\end{predl}

{\bf Proof.} Since ${\bf K}'(k)>0$ for any $k$ it follows from
(\ref{2.30}) that for $q\to 1-0$  ${\bf K}(k)\to\infty$, and hence
$\lim_{q\to 1-0}k=1-0.$ In this case ${\bf K}'(1)=\pi/2$ and 
$\dn u=1$. Now from (\ref{2.31}) we have
$$
\lim_{q\to 1-0}{\bf Q}_\nu=-\lim_{q\to 1-0}\frac{1-q}{\ln q}
{\bf K}'(k)\dn[{\frac{2\ln q^{-\nu+1/2}}{\pi}}{\bf K}'(k)]=
\frac{\pi}{2}.
$$

\section{Some properties of the $q$-binomial formula}
\setcounter{equation}{0}

There is a $q$-analog of the classical binomial formula \cite{GR}
$$
(1-z)^{-a}=\sum_{k=0}^\infty\frac{(a)_k}{k!}z^k,~(a)_k=
\frac{\Gamma(a+k)}{\Gamma(a)}),
\qquad |z|<1,
$$
$$
\frac{(q^\alpha z,q)_\infty}{(z,q)_\infty}=\sum_{k=0}^\infty
\frac{(q^\alpha,q)_k}{(q,q)_k}z^k,\qquad |z|<1.
$$
We need in two generalizations of the $q$-binom
\begin{equation}\label{3.1}
r(a,b,z,q)=\frac{(az,q)_\infty}{(bz,q)_\infty}
\end{equation}
\begin{equation}\label{3.2}
R(a,b,\gamma,z,q^2)=\frac{(az^2,q^2)_\infty}
{(bz^2,q^2)_\infty}z^\gamma
\end{equation}
\begin{predl}\label{p3.1}
The function $r(a,b,z,q)$ (\ref{3.1}) satisfies the difference equation
\begin{equation}\label{3.3}
z[br(a,b,z,q)-ar(a,b,qz,q)]=r(a,b,z,q)-r(a,b,qz,q).
\end{equation}
\end{predl}

{\bf Proof.} It follows from (\ref{3.3}) that
$$
\frac{r(a,b,z,q)}{r(a,b,qz,q)}=\frac{1-az}{1-bz}.
$$
Now the statement of the Proposition follows from
(\ref{2.7}).\rule{5pt}{5pt}
\bigskip
\begin{predl}\label{p3.2}
The function $R(a,b,\gamma,z,q^2)$ (\ref{3.2}) satisfies
the difference equation
\begin{equation}\label{3.4}
z^2[bq^\gamma R(a,b,\gamma,z,q^2)-aR(a,b,\gamma,qz,q^2)]=
q^\gamma R(a,b,\gamma,z,q^2)-R(a,b,\gamma,qz,q^2).
\end{equation}
\end{predl}

{\bf Proof.} Take $z^{-\gamma}R(a,b,\gamma,z,q)=r(a,b,z^2,q^2)$,
and substitute it
in (\ref{3.3}). Then the equation
$$
z^2[bz^{-\gamma}R(a,b,\gamma,z,q^2)-
aq^{-\gamma}z^{-\gamma}R(a,b,\gamma,qz,q^2)]=
z^{-\gamma}R(a,b,\gamma,z,q)-q^{-\gamma}z^{-\gamma}R(a,b,\gamma,qz,q^2).
$$
leads directly to (\ref{3.4}).\rule{5pt}{5pt}
\bigskip
\begin{lem}
\label{l3.1}
If $|a|<|b|$ the function $r(a,b,z,q)$ can be represented
as the sum of the
partial functions
\begin{equation}
\label{3.5}
\frac{(az,q)_\infty}{(bz,q)_\infty}=\frac1{(q,q)_\infty}
\sum_{k=0}^\infty\frac{(-1)^kq^{\frac{k(k+1)}2}(a/bq^{-k},q)_\infty}
{(q,q)_k(1-zbq^k)}.
\end{equation}
The series (\ref{3.5}) converges absolutely for any $z\ne
b^{-1}q^{-k},\quad k=0, 1,\ldots$
\end{lem}

{\bf Proof.} Let
$$
r_n(a,b,z,q)=\frac{(az,q)_n}{(bz,q)_n}=(a/b)^n+
\sum_{k=0}^n\frac{c_{k,n}}{1-zbq^k}.
$$
Then
$$
c_{k,n}=\lim_{z\to b^{-1}q^{-k}}(1-zbq^k)r_n(a,b,z,q)=
\frac{(-1)^kq^{\frac{k(k+1)}2}(a/bq^{-k},q)_n}
{(q,q)_k(q,q)_{n-k}}
$$
and
$$
r_n(a,b,z,q)=(a/b)^n+
\sum_{k=0}^n\frac{(-1)^kq^{\frac{k(k+1)}2}(a/bq^{-k},q)_n}
{(q,q)_k(1-zbq^k)}\frac1{(q,q)_{n-k}}.
$$
Since $\lim_{n\to\infty}(a/b)^n=0$ and
$\frac1{(q,q)_{n-k}}<\frac1{(q,q)_\infty}$ we obtain
$$
r(a,b,z,q)=\lim_{n\to\infty}r_n(a,b,z,q)=\frac1{(q,q)_\infty}
\sum_{k=0}^\infty\frac{(-1)^kq^{\frac{k(k+1)}2}(a/bq^{-k},q)_\infty}
{(q,q)_k(1-zbq^k)}.
$$
Applying the d`Alembert criterion  we obtain
$$
\lim_{k\to\infty}|\frac{q^{\frac{(k+1)(k+2)}2}(a/bq^{-k-1},q)_\infty
(q,q)_k(1-zbq^k)}
{(q,q)_{k+1}(1-zbq^{k+1})q^{\frac{k(k+1)}2}(a/bq^{-k},q)_\infty}|=
\lim_{k\to\infty}|q^{k+1}-a/b|=|a/b|<1.
$$
Hence (\ref{3.5}) converges absolutely for any $|a|<|b|$, and 
$r(a,b,z,q)$
is the meromorphic function with the ordinary poles at the points
$z=b^{-1}q^{-k},\quad k=0, 1,\ldots$ \rule{5pt}{5pt}
\bigskip

It is easy to get
\begin{rem}
\label{r3.1}
If $0<|a|<|b|$ then
\begin{equation}\label{3.6}
\frac{(az,q)_\infty}{(bz,q)_\infty}=
\frac{(a/b,q)_\infty}{(q,q)_\infty}
\sum_{k=0}^\infty\frac{(b/aq,q)_k(a/b)^k}{(q,q)_k(1-zbq^k)}.
\end{equation}
If $a=0$ then
\begin{equation}\label{3.6a}
\frac1{(bz,q)_\infty}=\frac1{(q,q)_\infty}
\sum_{k=0}^\infty\frac{(-1)^kq^{\frac{k(k+1)}2}}
{(q,q)_k(1-zbq^k)}.
\end{equation}
\end{rem}
\bigskip
Assume that $a=\epsilon q^{2\alpha}, b=\epsilon q^{2\beta},
\epsilon=\pm1$ in
(\ref{3.2}).  Then we have from (\ref{3.6})
\begin{cor}
\label{c3.1}
\begin{equation}\label{3.7}
z^\gamma\frac{(\epsilon q^{2\alpha}z^2,q^2)_\infty}
{(\epsilon q^{2\beta}z^2,q^2)_\infty}=
z^\gamma\frac{(q^{2(\alpha-\beta)},q^2)_\infty}
{(q^2,q^2)_\infty}\sum_{k=0}^\infty\frac
{(q^{2(\beta-\alpha+1)},q^2)_kq^{2(\alpha-\beta)k}}
{(q^2,q^2)_k(1-\epsilon z^2q^{2(\beta+k)})}.
\end{equation}
\end{cor}
\begin{rem}
\label{r3.2}
As it follows from \cite{GR}(1.3.2)
\begin{equation}\label{3.8}
\frac{(\epsilon q^{2\alpha}z^2,q^2)_\infty}
{(\epsilon q^{2\beta}z^2,q^2)_\infty}=\sum_{k=0}^\infty
\epsilon^kq^{2\beta k}\frac{(q^{2(\alpha-\beta)},q^2)_k}
{(q^2,q^2)_k}z^{2k},
\end{equation}
which converges in the domain $|z|<q^{-\beta/2}$.
\end{rem}

It follows from Lemma \ref{l3.1} that if $\gamma=0$ (\ref{3.7}) 
is the
meromorphic function with the ordinary poles $z=\pm\sqrt{\epsilon}
q^{-\beta-k},\quad k=0, 1,\ldots$, and hence it is the analitic
continuation of (\ref{3.8}).

\begin{cor}
\label{c3.2}
For an arbitrary real $s\ne0$
$$
\lim_{m\to\infty}|e_q(i\qpop q^{-m}s)|=0.
$$
\end{cor}

{\bf Proof.} Due to (\ref{2.9a}),  Lemma \ref{l3.1} for $a=0$,
$b=~i\qpop q^{-m}$,
and the evident inequality for $k\geq 0$
\begin{equation}\label{3.9}
\frac1{1+s^2(\qpop)^2q^{-2m+2k}}\le\frac{q^{-2k}}{1+s^2(\qpop)^2q^{-2m}},
\end{equation}
we have
$$
|e_q(i\qpop q^{-m}s)|\le\frac1{(q,q)_\infty}\sum_{k=0}^\infty
|\frac{q^{\frac{k(k+1)}2}}{(q,q)_k(1-i\qpop sq^{-m+k})}|\le
$$
$$
\le(1+s^2(\qpop)^2q^{-2m})^{-1/2}\frac1{(q,q)_\infty}\sum_{k=0}^\infty
\frac{q^{\frac{k(k-1)}2}q^{-k}}{(q,q)_k}=(1+s^2(\qpop)^2q^{-2m})^{-1/2}
E_q(q^{-1})e_q(q).
$$\rule{5pt}{5pt}
\bigskip

There are two types of  $q$-trigonometric functions
$$
\cos_qz=\frac12[e_q(iz)+e_q(-iz)],
$$
$$
\Cos_qz=\frac12[E_q(iz)+E_q(-iz)],\quad
\Sin_qz=\frac1{2i}[E_q(iz)-E_q(-iz)].
$$
\begin{cor}\label{c3.3}
For  real $s\ne0$ and  integer $m$
$$
|\cos_q(\qpop q^{-m}s)|\le\frac{E_q(q^{-1})e_q(q)}
{1+(\qpop)^2q^{-2m}s^2}.
$$
\end{cor}

{\bf Proof.} If $a=0,\quad
b=\pm i\qpop q^{-m}$ then from (\ref{3.6a})
$$
\cos_q(\qpop q^{-m}s)=\frac1{(q,q)_\infty}\sum_{k=0}^\infty
\frac{(-1)^kq^{\frac{k(k+1)}2}}{(q,q)_k(1+(\qpop)^2q^{-2m+2k}s^2)}.
$$

Using  (\ref{3.9}) we can put
$$
|\cos_(\qpop q^{-m}s)|\le\frac1{(q,q)_\infty(1+(\qpop)^2q^{-2m}s^2)}
\sum_{k=0}^\infty\frac{q^{\frac{k(k+1)}2}q^{-2k}}{(q,q)_k}.
$$
Then the statement follows from (\ref{2.9b})
\rule{5pt}{5pt}
\bigskip
\begin{cor}\label{c3.4}
For an arbitrary real $s>0$
$$
\lim_{m\to\infty}e_q(-\qpop q^{-m}s)=0.
$$
\end{cor}

{\bf Proof.} Again, due to (\ref{2.9a}), (\ref{2.9b}),
Lemma \ref{l3.1} for $a=0$, $b=~-\qpop q^{-m}$,
and the inequality
$$
\frac1{1+s\qpop q^{-m+k}}\le\frac{q^{-k}}{1+s\qpop q^{-m}},
$$
we have
$$
|e_q(-\qpop q^{-m}s)|\le\frac1{(q,q)_\infty}\sum_{k=0}^\infty
|\frac{q^{\frac{k(k+1)}2}}{(q,q)_k(1+\qpop sq^{-m+k})}|\le
$$
$$
\le\frac1{|1+s\qpop q^{-m}|}\frac1{(q,q)_\infty}\sum_{k=0}^\infty
\frac{q^{\frac{k(k-1)}2}}{(q,q)_k}=\frac1{|1+s\qpop q^{-m}|}
E_q(1)e_q(q).\rule{5pt}{5pt}
$$
\bigskip
\begin{cor}\label{c3.5}
For an arbitrary real $s$
\begin{equation}\label{3.10}
|\Cos_q(\qpop q^{-m}s)|\le 1,\quad
|\Sin_q(\qpop q^{-m}s)|\le\frac12q^{-m}(1+q)|s|.
\end{equation}
\end{cor}

{\bf Proof.} These functions are represented by the absolutely
convergent alternating series
$$
\Cos_q(\qpop q^{-m}s)=\sum_{n=0}^\infty
\frac{(-1)^nq^{n(2n-1-2m)}(1-q^2)^{2n}s^{2n}}{2^{2n}(q,q)_{2n}},
$$
$$
\Sin_q(\qpop q^{-m}s)=\sum_{n=0}^\infty
\frac{(-1)^nq^{(n-m)(2n+1)}(1-q^2)^{2n+1}s^{2n+1}}
{2^{2n+1}(q,q)_{2n+1}},
$$
and hence we have (\ref{3.10}).\rule{5pt}{5pt}
\bigskip
\begin{cor}\label{c3.6}
If $\alpha>\beta+1$ and real
$z\ne0$, then
$$
\frac{(-q^{2\alpha}z^2,q^2)_\infty}{(-q^{2\beta}z^2,q^2)_\infty}\le
\frac{C_{\alpha,\beta}}{1+z^2q^{2\beta}}.
$$
\end{cor}

{\bf Proof.} Since for any $k\geq 0\qquad \frac1{1+z^2q^{2\beta+2k}}\le
\frac{q^{-2k}}{1+z^2q^{2\beta}}$  we have from (\ref{3.7})
$$
\frac{(-q^{2\alpha}z^2,q^2)_\infty}{(-q^{2\beta}z^2,q^2)_\infty}\le
\frac{(q^{2(\alpha-\beta)},q^2)_\infty}
{(q^2,q^2)_\infty(1+z^2q^{2\beta})}\sum_{k=0}^\infty\Bigl|\frac
{(q^{2(\beta-\alpha+1)},q^2)_kq^{2(\alpha-\beta-1)k}}
{(q^2,q^2)_k}\Bigr|.
$$
The series in the right hand side  converges due
to the d'Alembert criterion
and thereby produce the finite contribution in the constant
 $C_{\alpha,\beta}$.
\rule{5pt}{5pt}
\bigskip
\begin{cor}\label{c3.7}
If $\alpha>\beta+1/2$
\begin{equation}\label{3.11}
\lim_{m\to\infty}\frac{(-q^{2(\alpha-m)},q^2)_\infty}
{(-q^{2(\beta-m)},q^2)_\infty}=0.
\end{equation}
\end{cor}

{\bf Proof.} For $z=q^{-2m}, \gamma=0,
\epsilon=-1$
 we have from (\ref{3.7}) $$
\frac{(-q^{2(\alpha-m)},q^2)_\infty}
{(-q^{2(\beta-m)},q^2)_\infty}=\frac{(q^{2(\alpha-\beta)},q^2)_\infty}
{(q^2,q^2)_\infty}\sum_{k=0}^\infty\frac{(q^{2(\beta-\alpha+1},q^2)_k
q^{2(\alpha-\beta)k}}{(q^2,q^2)_k(1+q^{2(\beta-m+k)})}.
$$
Obviously
$$
\frac1{1+q^{2(\beta-m+k)}}=\frac{q^{m-k-\beta}}{q^{m-k-\beta}+
q^{\beta+k-m}}\le\frac{q^{m-k-\beta}}2.
$$
Thus
$$
|\frac{(-q^{2(\alpha-m)},q^2)_\infty}
{(-q^{2(\beta-m)},q^2)_\infty}|\le\frac{q^{m-\beta}}2
\sum_{k=0}^\infty|\frac{(q^{2(\beta-\alpha+1)},q^2)_k
q^{2(\alpha-\beta-1/2)k}}{(q^2,q^2)_k}|.
$$
The series in the right hand side  converges. Due to the
multiplier $q^m$ we come to (\ref{3.11}).\rule{5pt}{5pt}
\bigskip
\begin{rem}\label{r3.3}
Let $a=\epsilon q^{2\alpha}, b=\epsilon q^{2\beta}$ in (\ref{3.2}) and
(\ref{3.4}). Then if $q\to1-0$ the difference equation (\ref{3.4})
takes the form of the differential equation
\begin{equation}\label{3.12}
z(1-\epsilon z^2)R'(z)-[\gamma+\epsilon(2\alpha-2\beta-\gamma)z^2]R(z)
=0
\end{equation}
with solution
$$
R(z)=Cz^\gamma(1-\epsilon z^2)^{\beta-\alpha}.
$$
\end{rem}

\section{The Jackson integral representation of the Modified $q$-Bessel
Functions}
\setcounter{equation}{0}

We remind that  modified $q$-Bessel function of the first kind
 \cite{OR1}
\begin{equation}\label{4.a}
\mbes=
\frac1{\G_{q^2}(\nu+1)}\sum_{n=0}^\infty\frac{(1-q^2)^{2n}s^{\nu+2n}}
{(q^2,q^2)_n(q^{2\nu+2},q^2)_n2^{\nu+2n}},\quad |s|<\frac1{1-q^2},
\end{equation}
satisfies the difference equation

\begin{equation}
\label{4.1}
[1-(\qpop)^2q^{-2}s^2]f_\nu(q^{-1}s)-(q^{-\nu}+q^\nu)f_\nu(s)+
f_\nu(qs)=0,
\end{equation}
and  modified $q$-Bessel function of the second kind
\begin{equation}\label{4.b}
\mbess=\frac1{\G_{q^2}(\nu+1)}\sum_{n=0}^\infty
\frac{q^{2n(\nu+n)}(1-q^2)^{2n}s^{\nu+2n}}
{(q^2,q^2)_n(q^{2\nu+2},q^2)_n2^{\nu+2n}}
\end{equation}
 satisfies the difference equation
\begin{equation}
\label{4.2}
f_\nu(q^{-1}s)-(q^{-\nu}+q^\nu)f_\nu(s)+
[1-(\qpop)^2s^2]f_\nu(qs)=0.
\end{equation}

$q$-Bessel-Macdonald function ($q$-BMF) has been determined for
noninteger $\nu$ as \cite{OR1}
\begin{equation}\label{5.1}
\makdj=\frac{q^{-\nu^2+1/2}}{4(a_\nu a_{-\nu})^{3/2}\sin\nu\pi}
[a_\nu I_{-\nu}^{(j)}((1-q^2)s;q^2)-a_{-\nu}I_\nu^{(j)}((1-q^2)s;q^2)]
\end{equation}
where
\begin{equation}\label{5.1a}
a_\nu=\sqrt{2/(1-q^2)}e_q(-1)\frac{I_\nu^{(2)}(2;q^2)}
{\phantom._2\Phi_1(q^{\nu+1/2},q^{-\nu+1/2};-q;q,q)},
\end{equation}
and
\begin{equation}\label{5.1b}
a_\nu a_{-\nu}=\frac{q^{-\nu+1/2}}
{2\G_{q^2}(\nu)\G_{q^2}(1-\nu)\sin\nu\pi}.
\end{equation}
If $\nu$ is an integer $n$  functions $K_n^{(j)}$ are well defined
 in
 the limit for $\nu\to n$ in (\ref{5.1}).

 Functions (\ref{5.1}) satisfy
 (\ref{4.1}) and (\ref{4.2}) for $j=1$ and $j=2$
respectively.

\begin{predl}
\label{p4.1}
Modified $q$-Bessel function ($q$-MBF) $I_\nu^{(1)}$ for $\nu>0$ can be
represented as the $q$-integral
\begin{equation}
\label{4.3}
\mbes=\frac{1+q}{2\G_{q^2}(\nu+1/2)\G_{q^2}(1/2)}
\int_{-1}^1d_qz\frac{(q^2z^2,q^2)_\infty}{(q^{2\nu+1}z^2,q^2)_\infty}
\Eex(s/2)^\nu.
\end{equation}
\end{predl}

{\bf Proof.} Consider the $q$-integral
\begin{equation}
\label{4.4}
S_1^{(1)}(s)=\int_{-1}^1d_qzf_\nu^{(1)}(z)\Eex,
\end{equation}
where $f_\nu^{(1)}(z)$ is a such function that it is absolutely convergent.
Require that $S_1^{(1)}(s)(s/2)^\nu$ satisfies (\ref{4.1}).  
Then $S_1^{(1)}(s)$ satisfies the equation
$$
[1-(\qpop)^2q^{-2}s^2]q^{-\nu}S_1^{(1)}(q^{-1}s)-
(q^{-\nu}+q^\nu)S_1^{(1)}(s)+q^\nu S_1^{(1)}(qs)=0.
$$
or
\begin{equation}
\label{4.5}
S_1^{(1)}(q^{-1}s)-S_1^{(1)}(s)-q^{2\nu}[S_1^{(1)}(s)-S_1^{(1)}(qs)]=
(\qpop)^2q^{-2}S_1^{(1)}(q^{-1}s)s^2.
\end{equation}
Substituting (\ref{4.4}) in (\ref{4.5}), multiplying it on
$\frac{2s^{-1}}{1-q^2}$ from the right, and drawing it
through $E_q$ we obtain
$$
\int_{-1}^1d_qzf_\nu^{(1)}(z)\frac{2s^{-1}}{1-q^2}[\Eexqq-\Eexq]-
$$
$$
-q^{2\nu}\int_{-1}^1d_qzf_\nu^{(1)}(z)\frac{2s^{-1}}{1-q^2}[\Eexq-\Eex]=
$$
$$
=\qpop q^{-2}\int_{-1}^1d_qzf_\nu^{(1)}(z)\Eexq s.
$$
Due to (\ref{2.11}), (\ref{2.18}) and (\ref{2.19}) it can be rewritten as
$$
q^{-1}\int_{-1}^1d_qzf_\nu^{(1)}(z)z\exoq-
q^{2\nu}\int_{-1}^1d_qzf_\nu^{(1)}(z)z\exo=
$$
$$
=\qpop q^{-2}\int_{-1}^1d_qzf_\nu^{(1)}(z)\exoq s,
$$
or
$$
\int_{-1}^1d_qzf_\nu^{(1)}(z)z^{-2\nu+1}\frac{2z^{-1}}{1-q^2}
[z^{2\nu+1}\exoq-q^{2\nu+1}z^{2\nu+1}\exo]=
$$
$$
=\int_{-1}^1d_qzf_\nu^{(1)}(z)\proz\exoq.
$$
Using  the $q$-integration by parts (\ref{2.2}) we obtain
$$
-\int_{-1}^1d_qz\p_z(f_\nu^{(1)}(z)z^{-2\nu+1})q^{2\nu+1}z^{2\nu+1}\exo=
-\int_{-1}^1d_qz\p_zf_\nu^{(1)}(z)\exo.
$$
Thus we come to the difference equation for $f_\nu^{(1)}(z)$
\begin{equation}
\label{4.6}
q^{2\nu+1}z^2[f_\nu^{(1)}(z)-q^{-2\nu+1}f_\nu^{(1)}(qz)]=
f_\nu^{(1)}(z)-f_\nu^{(1)}(qz).
\end{equation}
It coincides with (\ref{3.4}) for $a=q^2, b=q^{2\nu+1},\gamma=0$, and
hence from the Proposition \ref{p3.2}
\begin{equation}
\label{4.7}
f_\nu^{(1)}(z)=\frac{(q^2z^2,q^2)_\infty}{(q^{2\nu+1}z^2,q^2)_\infty}.
\end{equation}
 $S_1^{(1)}(s)(s/2)^\nu$ is a solution to  (\ref{4.1}) and therefore
 it can be represented as \cite{OR1}
$$
S_1^{(1)}(s)(s/2)^\nu=A\mbes+B\makd.
$$
Multiplying the both sides on $(s/2)^\nu$ and putting $s=0$ from
(\ref{4.a}) and (\ref{5.1}) we
obtain $B=0$. Multiplying again on $(s/2)^{-\nu}$ and assuming
$s=0$ we come to
$$
\int_{-1}^1d_qzf_\nu^{(1)}(z)=A\frac1{\G_{q^2}(\nu+1)}.
$$
To calculate the $q$-integral in the left hand side we use the
obvious property
$$
\int_{-1}^1f(x^2)d_qx=2(1-q)\sum_{m=0}^\infty q^mf(q^{2m})=
\frac{2(1-q^2)}{1+q}\sum_{m=0}^\infty q^{2m}f(q^{2m})q^{-m}=
$$
$$
=\frac2{1+q}\int_0^1f(x)x^{-1/2}d_{q^2}x.
$$
Then from \cite{GR} (1.11.7)
$$
\int_{-1}^1\frac{(q^2z^2,q^2)_\infty}{(q^{2\nu+1}z^2,q^2)_\infty}d_qz=
\frac2{1+q}\int_0^1\frac{(q^2z,q^2)_\infty}{(q^{2\nu+1}z,q^2)_\infty}
z^{-1/2}d_{q^2}z=\frac2{1+q}B_{q^2}(\nu+1/2,1/2),
$$
where $B_{q}(\nu;\mu)=\frac{\Gamma_q(\nu)\Gamma_q(\mu)}{\Gamma_q(\nu+\mu)}$
is  the $q$-beta function.
Hence
\begin{equation}\label{4.8}
A=\frac2{1+q}B_{q^2}(\nu+1/2,1/2)\G_{q^2}(\nu+1)=
\frac2{1+q}\G_{q^2}(\nu+1/2)\G_{q^2}(1/2).
\end{equation}
and we come to (\ref{4.3}).\rule{5pt}{5pt}
\bigskip
\begin{predl}
\label{p4.2}
The $q$-MBF $\mbess$ has the following $q$-integral representation
$$
\mbess=\frac1{2\G_{q^2}(\nu+1)
\phantom._2\Phi_1(q^{-2\nu+1},q;q^3;q^2,1)}\times
$$
\begin{equation}\label{4.9}
\times\int_{-1}^1d_qz\frac{(q^{-2\nu+1}z^2,q^2)_\infty}
{(z^2,q^2)_\infty} \bg(s/2)^\nu.
\end{equation}
\end{predl}

{\bf Proof.} Consider as above the absolutely convergent $q$-integral
\begin{equation}\label{4.10}
S_1^{(2)}(s)=\int_{-1}^1d_qzf_\nu^{(2)}(z)\bg.
\end{equation}
and assume that $S_1^{(2)}(s)(s/2)^\nu$ satisfies
(\ref{4.2}). Then $S_1^{(2)}(s)$ satisfies the equation
$$
q^{-\nu}S_1^{(2)}(q^{-1}s)-(q^{-\nu}+q^\nu)S_1^{(2)}(s)+
[1-(\qpop)^2s^2]q^\nu S_1^{(2)}(qs)=0,
$$
or
$$
q^{-2\nu}[S_1^{(2)}(q^{-1}s)-S_1^{(2)}(s)]-
[S_1^{(2)}(s)-S_1^{(2)}(qs)]=(\qpop)^2S_1^{(2)}(qs)s^2.
$$

The further arguments are the same as  in the proof of
Preposition \ref{p4.1}. Here we use
(\ref{2.10}),(\ref{2.20}) and (\ref{2.21}) and arrive to the difference
equation for $f_\nu^{(2)}(z)$:
\begin{equation}\label{4.11}
z^2[f_\nu^{(2)}(z)-q^{-2\nu+1}f_\nu^{(2)}(qz)]=
f_\nu^{(2)}(z)-f_\nu^{(2)}(qz).
\end{equation}
It coincides with (\ref{3.4}) for $a=q^{-2\nu+1}, b=1, \gamma=0$,
and hence from  Proposition \ref{p3.2}
$$
f_\nu^{(2)}(z)=\frac{(q^{-2\nu+1}z^2,q^2)_\infty}{(z^2,q^2)_\infty}.
$$
Since $S_1^{(2)}(s)(s/2)^\nu$ is a solution of (\ref{4.2})
it can be represented as
$$
S_1^{(2)}(s)(s/2)^\nu=A\mbess+B\makdo.
$$
 As in the proof of Preposition \ref{p4.1} from (\ref{4.b}) and
(\ref{5.1}) we obtain $B=0$ and
$$
A=\G_{q^2}(\nu+1)\int_{-1}^1\frac{(q^{-2\nu+1}z^2,q^2)_\infty}
{(z^2,q^2)_\infty}d_qz=\frac{2\G_{q^2}(\nu+1)}{1+q}
\int_0^1\frac{(q^{-2\nu+1}z,q^2)_\infty}{(z,q^2)_\infty}z^{-1/2}d_{q^2}z.
$$
It follows from \cite{GR} (1.11.9) that
$$
A=2\G_{q^2}(\nu+1)\phantom._2\Phi_1(q^{-2\nu+1},q;q^3;q^2,1),
$$
and we come to (\ref{4.9}).\rule{5pt}{5pt}
\bigskip
\begin{rem}\label{r4.1}
 If $q\to1-0$ the equations
(\ref{4.6}) and (\ref{4.11}) take the form of the differential equation
(see  Remark \ref{r3.3})
$$
(1-z^2)f_\nu'(z)+(2\nu-1)zf_\nu(z)=0.
$$
The solution to this equation is
$$
f_\nu(z)=C(1-z^2)^{\nu-1/2},
$$
which leads to the classical integral representation of  Modified
Bessel function \cite{BE} (7.12.10)
$$
I_\nu(s)=\frac{(s/2)^\nu}{\G(\nu+1/2)\G(1/2)}
\int_{-1}^1(1-z^2)^{\nu-1/2}e^{zs}dz.
$$
\end{rem}
\bigskip
As it follows from
$$
I_\nu^{(j)}((1-q^2)s;q^2)=
e^{-i\nu\pi/2}J_\nu^{(j)}((1-q^2)e^{i\pi/2}s;q^2),\quad j=1, 2.
$$
Then (\ref{4.3}) and (\ref{4.9}) give the $q$-integral representations
\begin{cor}\label{c4.1}
$$
J_\nu^{(1)}((1-q^2)s;q^2)=\frac{(1+q)}
{2\Gamma_{q^2}(\nu+1/2)\Gamma_{q^2}(1/2)}\int_{-1}^1d_qz
\frac{(q^2z^2,q^2)_\infty}{(q^{2\nu+1}z^2,q^2)_\infty}
E_q(-i\qpop zs)(s/2)^\nu,
$$
$$
J_\nu^{(2)}((1-q^2)s;q^2)=\frac{(1+q)}
{2\Gamma_{q^2}(\nu+1)\phantom._2\Phi_1(q^{-2\nu+1},q;q^3;q^2,1)}
\times
$$
\begin{equation}\label{4.12}
\times\int_{-1}^1d_qz
\frac{(q^{-2\nu+1}z^2,q^2)_\infty}{(z^2,q^2)_\infty}
\phantom._0\Phi_1(-;0;q,-i\qpop zs)(s/2)^\nu.
\end{equation}
\end{cor}

\section{The $q$-Fourier integral representation of
$q$-Bessel-Macdonald Functions}
\setcounter{equation}{0}

Here we consider $q$-BMF (\ref{5.1}) . We will use the following
formulas from \cite{OR1}:
\begin{equation}\label{5.a}
\mbess=\frac{a_\nu}{\sqrt{z}}[E_q(\qpop s)\Phi_\nu(s)+ie^{i\nu\pi}
E_q(-\qpop s)\Phi_\nu(-s)],
\end{equation}
\begin{equation}\label{5.b}
\makdo=\frac{q^{-\nu^2+1/2}}{2\sqrt{a_\nu a_{-\nu}s}}
E_q(-\qpop s)\Phi_\nu(s),
\end{equation}
where
$$
\Phi_\nu(s)=\phantom._2\Phi_1(q^{\nu+1/2},q^{-\nu+1/2};-q;q,
\frac{2q}{(1-q^2)s}).
$$

\begin{predl}\label{p5.1}
 $q$-BMF $\makd$ for $\nu>1/2$ can be represented by the $q$-integral
$$
\makd=\frac{q^{-\nu^2+1/2}\Gamma_{q^2}(\nu+1/2)\Gamma_{q^2}(1/2)}
{4{\bf Q}_\nu}\sqrt{\frac{a_\nu}{a_{-\nu}}}
\times
$$
\begin{equation}\label{5.2}
\times\int_{-\infty}^\infty d_qz\frac{(-q^2z^2,q^2)_\infty}
{(-q^{-2\nu+1}z^2,q^2)_\infty}E_q(i\qpop zs)(s/2)^{-\nu},
\end{equation}
where ${\bf Q}_\nu$ is defined by (\ref{2.31}).
\end{predl}

{\bf Proof.} Consider the $q$-integral
\begin{equation}\label{5.3}
S_2^{(1)}(s)=\int_{-\infty}^\infty d_qzh_\nu^{(1)}(z)\Eexi,
\end{equation}
and assume that it absolutely converges together with its
$q$-derivative. According to  Definition \ref{d2.1} and
(\ref{2.11}) it means that
$$
\lim_{m\to\pm\infty}q^m|h_\nu^{(1)}(q^m)\frac{2\p_s}{1+q}
e_q(i\qpop q^{m-1}s)|=0.
$$
It follows from (\ref{2.19})
\begin{equation}\label{5.5a}
\lim_{m\to\pm\infty}q^{2m}|h_\nu^{(1)}(q^m)
e_q(i\qpop q^ms)|=0.
\end{equation}
Substitute $S_2^{(2)}(s)(s/2)^{-\nu}$
in (\ref{4.1}). Then $S_2^{(1)}(s)$ satisfies the equation
\begin{equation}\label{5.4}
S_2^{(1)}(q^{-1}s)-S_2^{(1)}(s)-q^{-2\nu}[S_2^{(1)}(s)-
S_2^{(1)}(qs)]=(\qpop)^2q^{-2}S_2^{(1)}(q^{-1}s)s^2.
\end{equation}
Substituting (\ref{5.3}) in (\ref{5.4}), multiplying it on
$\frac{2s^{-1}}{1-q^2}$ from the right and drawing the multiplier
through $E_q$ we obtain
$$
\int_{-\infty}^\infty d_qzh_\nu^{(1)}(z)\frac{2s^{-1}}{1-q^2}
[\Eexiqq-\Eexiq]-
$$
$$
-\int_{-\infty}^\infty d_qzh_\nu^{(1)}(z)\frac{2s^{-1}}{1-q^2}
[\Eexiq-\Eexi]=
$$
$$
=\qpop q^{-2}\int_{-\infty}^\infty d_qzh_\nu^{(1)}(z)
\Eexiq s.
$$
Using (\ref{2.11}), (\ref{2.18}), and (\ref{2.19}) we come to
$$
\int_{-\infty}^\infty d_qzh_\nu^{(1)}(z)z[\exoiq-q^{-2\nu+1}\exoi]=
$$
$$
=-i\qpop q^{-1}\int_{-\infty}^\infty d_qzh_\nu^{(1)}(z)\exoiq s,
$$
or
$$
\int_{-\infty}^\infty
d_qzh_\nu^{(1)}(z)z^{2\nu+1}\frac{2z^{-1}}{1-q^2}[z^{-2\nu+1}\exoiq-
q^{-2\nu+1}z^{-2\nu+1}\exoi]=
$$
$$
=\int_{-\infty}^\infty d_qzh_\nu^{(1)}(z)\proz\exoiq.
$$
 $q$-Integration by parts (\ref{2.4}) gives
$$
\lim_{m\to\infty}[h_\nu^{(1)}(q^{-m})(q^{-2m}+1)e_q(i\qpop q^{-m-1}s)+
h_\nu^{(1)}(-q^{-m})(q^{-2m}+1)e_q(-i\qpop q^{-m-1}s)]=
$$
\begin{equation}\label{5.5}
=\int_{-\infty}^\infty d_qz[\p_z(h_\nu^{(1)}(z)z^{2\nu+1})
q^{-2\nu+1}z^{-2\nu+1}+\p_zh_\nu^{(1)}(z)]\exoi.
\end{equation}
According to (\ref{5.5a}) the left hand side of (\ref{5.5}) vanishes.
 Thus we come to the difference equation for $h_\nu^{(1)}(z)$
\begin{equation}\label{5.6}
z^2[-q^{-2\nu+1}h_\nu^{(1)}(z)+q^2h_\nu^{(1)}(qz)]=
h_\nu^{(1)}(z)-h_\nu^{(1)}(qz).
\end{equation}
Obviously
(\ref{5.6}) coincides with (\ref{3.4}) for $a=-q^2,
b=-q^{-2\nu+1}, \gamma=0$, and hence
\begin{equation}\label{5.7}
h_\nu^{(1)}(z)=\frac{(-q^2z^2,q^2)_\infty}{(-q^{-2\nu+1}z^2,q^2)_\infty}.
\end{equation}
It follows from  Corollaris \ref{3.2} , \ref{3.6} that
$h_\nu^{(1)}(z)$ satisfies (\ref{5.5a}). Then
$$
S_2^{(1)}(z)=\int_{-\infty}^\infty d_qz
\frac{(-q^2z^2,q^2)_\infty}{(-q^{-2\nu+1}z^2,q^2)_\infty}\Eexi.
$$
Since $S_2^{(1)}(s)(s/2)^{-\nu}$ is a solution of (\ref{4.1})
it can be represented in the form
$$
\int_{-\infty}^\infty d_qz
\frac{(-q^2z^2,q^2)_\infty}{(-q^{-2\nu+1}z^2,q^2)_\infty}\Eexi(s/2)^{-\nu}=
A\mbes+B\makd,
$$
or
$$
(1-q)\sum_{m=-\infty}^\infty q^m\frac{(-q^{2m+2},q^2)_\infty}
{(-q^{-2\nu+2m+1},q^2)_\infty}\cos_q(\qpop q^ms)(s/2)^{-\nu}=
$$
\begin{equation}\label{5.8}
=A\mbes+B\makd.
\end{equation}
As it follows from Corollary \ref{c3.3},
the function in the left hand side is holomorphic  in the domain
$\R s>0$, but $\mbes$
has  ordinary poles in the points $s=~\pm\frac{2q^{-r}}{1-q^2},\quad
r=0, 1,\ldots$ \cite{OR1}. Hence $A=0$. Multiplying (\ref{5.8}) on
$(s/2)^\nu$ and assuming $s=0$ we find that
$$
B=8a_{-\nu}\sqrt{a_\nu a_{-\nu}}q^{\nu^2-1/2}\sin\nu\pi\G_{q^2}(-\nu+1)
(1-q)\sum_{m=-\infty}^\infty q^m\frac{(-q^{2m+2},q^2)_\infty}
{(-q^{-2\nu+2m+1},q^2)_\infty}.
$$
Let calculate the last sum. It follows from the Corollary \ref{c3.1} for
$\alpha=m+1, \beta=-\nu+m+1/2,\\ \epsilon=-1, \gamma=0$
$$
\frac{(-q^{2m+2},q^2)_\infty}{(-q^{-2\nu+2m+1},q^2)_\infty}=
\frac{(q^{2\nu+1},q^2)_\infty}{(q^2,q^2)_\infty}
\sum_{k=0}^\infty\frac{q^{(2\nu+1)k}(q^{-2\nu+1},q^2)_k}
{(q^2,q^2)_k(1+q^{-2\nu+2m+2k+1})},
$$
and this series converges uniformly on real axis.

 So
$$
\sum_{m=-\infty}^\infty q^m\frac{(-q^{2m+2},q^2)_\infty}
{(-q^{-2\nu+2m+1},q^2)_\infty}=
\frac{(q^{2\nu+1},q^2)_\infty}{(q^2,q^2)_\infty}
\sum_{m=-\infty}^\infty q^m\sum_{k=0}^\infty
\frac{q^{(2\nu+1)k}(q^{-2\nu+1},q^2)_k}
{(q^2,q^2)_k(1+q^{-2\nu+2m+2k+1})}
$$
$$
=\frac{q^{\nu-1/2}(q^{2\nu+1},q^2)_\infty}{(q^2,q^2)_\infty}
\sum_{k=0}^\infty\frac{q^{2\nu k}(q^{-2\nu+1},q^2)_k}{(q^2,q^2)_k}
\sum_{m=-\infty}^\infty\frac1{q^{\nu-m-1/2}+q^{-\nu+m+1/2}}.
$$
Then using (\ref{2.7}), (\ref{2.31}) and (\ref{5.1a}), (\ref{5.1b}) we
obtain
$$
B=\frac{4q^{\nu^2+1/2}{\bf Q}_\nu}
{\Gamma_{q^2}(\nu+1/2)\Gamma_{q^2}(1/2)}
\sqrt{\frac{a_{-\nu}}{a_\nu}},
$$
and (\ref{5.2}).\rule{5pt}{5pt}
\bigskip
\begin{predl}\label{p5.2}
$q$-BMF $\makdo$ for $\nu>3/2$ can be represented by the $q$-integral
$$
\makdo=\frac{q^{-\nu^2+\nu}\Gamma_{q^2}(\nu+1/2)\Gamma_{q^2}(1/2)}
{4{\bf Q}_{1/2}}
\sqrt{\frac{a_\nu}{a_{-\nu}}}\times
$$
\begin{equation}\label{5.9}
\times\int_{-\infty}^\infty d_qz\frac{(-q^{2\nu+1}z^2,q^2)_\infty}
{(-z^2,q^2)_\infty}\ibg(s/2)^{-\nu},
\end{equation}
where ${\bf Q}_{1/2}$ is determined by (\ref{2.31}).
\end{predl}

{\bf Proof.} Consider the $q$-integral
\begin{equation}\label{5.10}
S_2^{(2)}(s)=\int_{-\infty}^\infty d_qzh_\nu^{(2)}(z)\ibg.
\end{equation}
Let it converges absolutely together with its
$q$-derivative. According to Definition \ref{d2.1} and
(\ref{2.10}) it means that
$$
\lim_{m\to\pm\infty}q^m|h_\nu^{(2)}(q^m)\frac{2\p_s}{1+q}
E_q(i\qpop q^{m-1}s)|=0
$$
It follows from (\ref{2.21})
\begin{equation}\label{5.12a}
\lim_{m\to\pm\infty}q^{2m}|h_\nu^{(2)}(q^m)
E_q(i\qpop q^ms)|=0.
\end{equation}

Require that $S_2^{(2)}(s)(s/2)^{-\nu}$ satisfies
(\ref{4.2}). Then $S_2^{(2)}(s)$ satisfies
\begin{equation}\label{5.10a}
q^{2\nu}[S_2^{(2)}(q^{-1}s)-S_2^{(2)}(s)]-[S_2^{(2)}(s)-
S_2^{(2)}(qs)]=(\qpop)^2S_2^{(2)}(qs)s^2.
\end{equation}
Using (\ref{2.10}), (\ref{2.20}) and (\ref{2.21}) and arguing as in
 the proof of
Proposition \ref{p5.1} we come to the equality
$$
q^{2\nu-1}\int_{-\infty}^\infty d_qzh_\nu^{(2)}(z)z^{2\nu+1}
\p_z(z^{-2\nu+1}\iExo)=-\int_{-\infty}^\infty
d_qzh_\nu^{(2)}(z)\partial_z\iExo.
$$
 $q$-Integration by parts (\ref{2.4}) gives
$$
\lim_{m\to\infty}q^{2\nu-1}[h_\nu^{(2)}(q^{-m})q^{-2m}
E_q(i\qpop q^{-m}s)+h_\nu^{(2)}(-q^{-m})q^{-2m}
E_q(-i\qpop q^{-m}s)]-
$$
$$
-\int_{-\infty}^\infty d_qz\p_z(h_\nu^{(2)}(z)z^{2\nu+1})z^{-2\nu+1}
\iExop=
$$
$$
=-\lim_{m\to\infty}[h_\nu^{(2)}(q^{-m})E_q(i\qpop q^{-m}s)
+h_\nu^{(2)}(-q^{-m}) E_q(-i\qpop q^{-m}s)]+
$$
$$
+\int_{-\infty}^\infty d_qz\p_zh_\nu^{(2)}(z)\iExop,
$$
or
$$
\lim_{m\to\infty}(q^{2\nu-1-2m}+1)[h_\nu^{(2)}(q^{-m})E_q(i\qpop q^{-m}s)
+h_\nu^{(2)}(-q^{-m}) E_q(-i\qpop q^{-m}s)]=
$$
\begin{equation}\label{5.11}
=\int_{-\infty}^\infty d_qz[\p_z(h_\nu^{(2)}(z)z^{2\nu+1})z^{-2\nu+1}+
\p_zh_\nu^{(2)}(z)]\iExop.
\end{equation}
(\ref{5.12a}) means that the left hand side of (\ref{5.11}) vanishes.
 Thus we come to the difference equation for $h_\nu^{(2)(z)}$
\begin{equation}\label{5.12}
z^2[h_\nu^{(2)}(z)-q^{2\nu+1}h_\nu^{(2)}(qz)]=
-h_\nu^{(2)}(z)+h_\nu^{(2)}(qz).
\end{equation}
It coincides with (\ref{3.4}) for $a=-q^{2\nu+1},
b=-1, \gamma=0$, and hence
\begin{equation}\label{5.13}
h_\nu^{(2)}(z)=\frac{(-q^{2\nu+1}z^2,q^2)_\infty}{(-z^2,q^2)_\infty}.
\end{equation}
It follows from Corollaris \ref{3.2} and \ref{3.6} that $h_\nu^{(2)}(z)$ 
satisfies 
(\ref{5.12a}). 
But $S_2^{(2)}(s)(s/2)^{-\nu}$ satisfies (\ref{4.2}) and thereby it is
represented in the form
$$
\int_{-\infty}^\infty d_qz\frac{(-q^{2\nu+1}z^2,q^2)_\infty}
{(-z^2,q^2)_\infty}\ibg(s/2)^{-\nu}=
$$
$$
=A\mbess+B\makdo,
$$
or
$$
2(1-q)\sum_{m=-\infty}^\infty q^m\frac{(-q^{2\nu+1+2m},q^2)_\infty}
{(-q^{2m},q^2)_\infty}\Cos_q(\qpop q^ms)(s/2)^{-\nu}=
$$
\begin{equation}\label{5.14}
=A\mbess+B\makdo.
\end{equation}

It follows from (\ref{5.a}) and (\ref{5.b}) that
        $\lim_{s\to\infty}\mbess=\infty$ and\\
$\lim_{s\to\infty}\makdo=0$.

Because the left hand side of (\ref{5.14}) vanishes
if $s\to\infty,\qquad A=0$. Multiplying (\ref{5.14}) on
$(s/2)^\nu$ and assuming $s=0$ we have
$$
B=8a_{-\nu}\sqrt{a_{-\nu}a_\nu}q^{\nu^2-1/2}\sin{\nu\pi}
\Gamma_{q^2}(-\nu+1)
(1-q)\sum_{m=-\infty}^\infty q^m\frac{(-q^{2\nu+1+2m},q^2)_\infty}
{(-q^{2m},q^2)_\infty}.
$$
Calculate the last sum. It follows from Corollary \ref{c3.1} for
$\alpha=\nu+m+1/2, \beta=m, \epsilon=-1, \gamma=0$ that
$$
\frac{(-q^{2\nu+1+2m},q^2)_\infty}{(-q^{2m},q^2)_\infty}=
\frac{(q^{2\nu+1},q^2)_\infty}{(q^2,q^2)_\infty}
\sum_{k=0}^\infty\frac{(q^{-2\nu+1},q^2)_kq^{(2\nu+1)k}}
{(q^2,q^2)_k(1+q^{2m+2k})}
$$
and this series converges uniformly on real axis.  So
$$
\sum_{m=-\infty}^\infty q^m\frac{(-q^{2\nu+1+2m},q^2)_\infty}
{(-q^{2m},q^2)_\infty}=
\frac{(q^{2\nu+1},q^2)_\infty}{(q^2,q^2)_\infty}
\sum_{m=-\infty}^\infty q^m
\sum_{k=0}^\infty\frac{(q^{-2\nu+1},q^2)_kq^{(2\nu+1)k}}
{(q^2,q^2)_k(1+q^{2m+2k})}=
$$
$$
=\frac{(q^{2\nu+1},q^2)_\infty}{(q^2,q^2)_\infty}
\sum_{k=0}^\infty\frac{(q^{-2\nu+1},q^2)_kq^{2\nu k}}
{(q^2,q^2)_k}\sum_{m=-\infty}^\infty\frac1{q^{-m}+q^m}
$$
Using (\ref{3.8}), (\ref{2.31}) and (\ref{5.1a}), (\ref{5.1b}) we obtain
$$
B=\frac{4q^{\nu^2-\nu}{\bf Q}_{1/2}}
{\G_{q^2}(\nu+1/2)\G_{q^2}(1/2)}
\sqrt{\frac{a_{-\nu}}{a_\nu}}
$$
and (\ref{5.9}).\rule{5pt}{5pt}
\bigskip
\begin{rem}\label{r5.1}
If $q\to 1-0$  equations (\ref{5.6}) and (\ref{5.12}) take the form
 of the
differential equation
$$
(z^2+1)h_\nu'(z)=-(2\nu+1)zh_\nu(z).
$$
with solution
$$
h_\nu(z)=C(z^2+1)^{-\nu-1/2},
$$
Using  Proposition \ref{p2.3} we obtain the classical integral
representation of  Bessel-Macdonald function (the Fourier integral)
\cite{BE} (7.12.27)
$$
K_\nu(s)=\frac{\G(\nu+1/2)(s/2)^{-\nu}}{2\G(1/2)}
\int_{-\infty}^\infty(z^2+1)^{-\nu-1/2}e^{izs}dz.
$$
\end{rem}

\section{ Bessel $q$-integral representation of
$q$-Bessel-Macdonald Functions}
\setcounter{equation}{0}
\begin{predl}\label{p6.1}
The $q$-BMF $\makd$ for $\nu>1/2$ can be represented as the $q$-integral
$$
\makd=\frac12q^{-\nu^2-\nu}(1+q)\Gamma_{q^2}(\nu+1)
\sqrt{\frac{a_\nu}{a_{-\nu}}}\times
$$
\begin{equation}\label{6.1}
\times\int_0^\infty d_qz\frac{(-q^2z^2,q^2)_\infty}
{(-q^{-2\nu}z^2,q^2)_\infty}z\iio(s/2)^{-\nu}.
\end{equation}
\end{predl}

{\bf Proof.} Consider the absolutely convergent $q$-integral
\begin{equation}\label{6.2}
S_3^{(1)}(s)=\int_0^\infty d_qzg_\nu^{(1)}(z)\iio
\end{equation}
  together with its
$q$-derivative. According to  Definition \ref{d2.1} and
(\ref{2.17}) it means
$$
\lim_{m\to\pm\infty}q^m|g_\nu^{(1)}(q^m)\frac{2\p_s}{1+q}
J_0^{(1)}((1-q^2)q^{m-1}s,q^2)|=0
$$
It follows from (\ref{2.23})
\begin{equation}\label{6.2a}
\lim_{m\to\pm\infty}q^{2m}|g_\nu^{(1)}(q^m)
J_1^{(1)}((1-q^2)q^ms)|=0.
\end{equation}

Assume that $S_3^{(1)}(s)(s/2)^{-\nu}$ satisfies
(\ref{4.1}). Then $S_3^{(1)}(s)$ satisfies  (\ref{5.4}).
Substituting (\ref{6.2}) in (\ref{5.4}), multiplying it on
$\frac{2s^{-1}}{1-q^2}$ from the right and drawing the multiplier
through $J_0^{(2)}$ on left we obtain
$$
\int_0^\infty d_qzg_\nu^{(1)}(z)\frac{2s^{-1}}{1-q^2}
[\iioqq-\iioq]-
$$
$$
-q^{-2\nu}\int_0^\infty d_qzg_\nu^{(1)}(z)\frac{2s^{-1}}{1-q^2}
[\iioq-\iio]=
$$
$$
=\qpop q^{-2}\int_0^\infty d_qzg_\nu^{(1)}(z)\iioq s.
$$
Using (\ref{2.17}), (\ref{2.22}) and (\ref{2.23}) we transform
it as
$$
-q^{-1}\int_0^\infty d_qzg_\nu^{(1)}(z)z\ioq+
q^{-2\nu}\int_0^\infty d_qzg_\nu^{(1)}(z)z\io=
$$
$$
=\qpop q^{-2}\int_0^\infty d_qzg_\nu^{(1)}(z)
\ddag J_0^{(1)}((1-q^2)q^{-1}zs;q^2)\ddag s,
$$
or
$$
\int_0^\infty d_qzg_\nu^{(1)}(z)z^{2\nu+1}\frac{2z^{-1}}{1-q^2}
[z^{-2\nu+1}\ioq-q^{-2\nu+1}z^{-2\nu+1}\io]=
$$
$$
=-q^{-1}\int_0^\infty d_qzg_\nu^{(1)}(z)
\ddag J_0^{(1)}((1-q^2)q^{-1}zs;q^2)\ddag s.
$$
It follows from (\ref{2.24}) that
$$
\int_0^\infty d_qzg_\nu^{(1)}(z)z^{2\nu+1}\p_z(z^{-2\nu+1}\ioq)=
$$
$$
=-\int_0^\infty d_qzg_\nu^{(1)}(z)(z/2)^{-1}\p_z(z/2\ioq).
$$
Using  the $q$-integration by parts (\ref{2.3}) we obtain
$$
\lim_{m\to\infty}[g_\nu^{(1)}
(q^{-m})q^{-2m}J_1^{(1)}((1-q^2)q^{-m-1}s;q^2)-
g_\nu^{(1)}(q^m)q^{2m}J_1^{(1)}((1-q^2)q^{m-1}s;q^2)]-
$$
$$
-\int_0^\infty d_qz\p_z(g_\nu^{(1)}(z)z^{2\nu+1})q^{-2\nu+1}z^{-2\nu+1}
\io=
$$
$$
=-\lim_{m\to\infty}[g_\nu^{(1)}(q^{-m})J_1^{(1)}((1-q^2)q^{-m-1}s;q^2)-
g_\nu^{(1)}(q^m)J_1^{(1)}((1-q^2)q^{m-1}s;q^2)]+
$$
$$
+\int_0^\infty d_qz\p_z(p_\nu^{(1)}(z)(z/2)^{-1})qz/2\io,
$$
or
$$
\lim_{m\to\infty}[g_\nu^{(1)}(q^{-m})(q^{-2m}+1)
J_1^{(1)}((1-q^2)q^{-m-1}s;q^2)-
g_\nu^{(1)}(q^m)(q^{2m}+1)J_1^{(1)}((1-q^2)q^{m-1}s;q^2)]=
$$
\begin{equation}\label{6.3}
=\int_0^\infty d_qz[\p_z(g_\nu^{(1)}(z)(z/2)^{-1})qz/2+
\p_z(g_\nu^{(1)}(z)z^{2\nu+1})q^{-2\nu+1}z^{-2\nu+1}]\io.
\end{equation}
It follows from (\ref{6.2a}) that the  left hand side of (\ref{6.3}) vanishes.
 Thus we come to the difference equation for $g_\nu^{(1)}(z)$
\begin{equation}\label{6.4}
z^2[q^{-2\nu+1}g_\nu^{(1)}(z)-q^2g_\nu^{(1)}(qz)]=
-qg_\nu^{(1)}(z)+g_\nu^{(1)}(qz).
\end{equation}
 But
(\ref{6.4}) coincides with (\ref{3.4}) for $a=-q^2, b=-q^{-2\nu},
\gamma=1$ and hence
\begin{equation}\label{6.5}
g_\nu^{(1)}(z)=\frac{(-q^2z^2,q^2)_\infty}{(-q^{-2\nu}z^2,q^2)_\infty}z.
\end{equation}

Since $S_3^{(1)}(s)(s/2)^{-\nu}$ satisfies  (\ref{4.1}) it can be
represented in the form
\begin{equation}\label{6.7}
\int_0^\infty d_qz\frac{(-q^2z^2,q^2)_\infty}
{(-q^{-2\nu}z^2,q^2)_\infty}z\iio(s/2)^{-\nu}=
A\mbes+B\makd.
\end{equation}
The function $\mbes$ is meromorphic  with  ordinary poles
$s=\pm\frac{2q^{-r}}{1-q^2}, r=0, 1,\ldots$. On the other hand the
function $\makd$ and the left hand side of (\ref{6.7})
are the holomorphic functions in the
domain $\R s>0$. Thereby $A=0$. Multiplying (\ref{6.7}) on $(s/2)^\nu$
and
taking $s=0$ we obtain
$$
\int_0^\infty
d_qz\frac{(-q^2z^2,q^2)_\infty}{(-q^{-2\nu}z^2,q^z)_\infty}z=
B\frac{q^{-\nu^2+1/2}}{4a_\nu\sqrt{a_\nu a_{-\nu}}\sin{\nu\pi}
\Gamma_{q^2}(-\nu+1)}.
$$
Rewriting the $q$-integral by means of (\ref{2.9a}) as
$$
\int_0^\infty d_qz\frac{(-q^2z^2,q^2)_\infty}
{(-q^{-2\nu}z^2,q^2)_\infty}z=\frac1{1+q}\int_0^\infty d_{q^2}z
\frac{(-q^2z,q^2)_\infty}{(-q^{-2\nu}z,q^2)_\infty}=
$$
$$
=\frac1{1+q}\int_0^\infty d_{q^2}zE_{q^2}(q^2z)e_{q^2}(-q^{-2\nu}z)
$$
and defining
$$
D_zf(z)=\frac{f(z)-f(q^2z)}{(1-q^2)z},
$$
$$
D_zE_{q^2}(z)=\frac1{1-q^2}E_{q^2}(q^2z),\quad
D_ze_{q^2}(-q^{-2\nu}z)=-\frac{q^{-2\nu}}{1-q^2}e_{q^2}(-q^{-2\nu}z).
$$
we obtain by means of (\ref{2.3}) 
$$
\int_0^\infty d_{q^2}zE_{q^2}(q^2z)e_{q^2}(-q^{-2\nu}z)=(1-q^2)
\int_0^\infty d_{q^2}zD_zE_{q^2}(z)e_{q^2}(-q^{-2\nu}z)=
$$
$$
=(1-q^2)\lim_{m\to\infty}[E_{q^2}(q^{-2m})e_{q^2}(-q^{-2\nu-2m})-
E_{q^2}(q^{2m})e_{q^2}(-q^{-2\nu+2m})]-
$$
$$
-(1-q^2)\int_0^\infty d_{q^2}zE_{q^2}(q^2z)D_ze_{q^2}(-q^{-2\nu}z)=
$$
$$
=(1-q^2)[\lim_{m\to\infty}\frac{(-q^{-2m},q^2)_\infty}
{(-q^{-2\nu-2m},q^2)_\infty}-1]+
q^{-2\nu}\int_0^\infty d_{q^2}zE_{q^2}(q^2z)e_{q^2}(-q^{-2\nu}z).
$$
Due to Corollary \ref{c3.7} for $\alpha=0, \beta=-\nu$ the limit vanishes
and we obtain
$$
\int_0^\infty d_z{q^2}zE_{q^2}(q^2z)e_{q^2}(-q^{-2\nu}z)=
-\frac{1-q^2}{1-q^{-2\nu}}.
$$
Hence
$$
\int_0^\infty d_qz\frac{(-q^2z^2,q^2)_\infty}
{(-q^{-2\nu}z^2,q^2)_\infty}z=-\frac{1-q}{1-q^{-2\nu}},
$$
and
$$
B=-\frac{1-q}{1-q^{-2\nu}}4q^{\nu^2-\frac12}a_{-\nu}
\sqrt{a_\nu a_{-\nu}}\sin{\nu\pi}\Gamma_{q^2}(-\nu+1)=
\frac{2q^{\nu^2+\nu}}{(1+q)\Gamma_{q^2}(\nu+1)}
\sqrt{\frac{a_\nu}{a_{-\nu}}}.
$$
Then from (\ref{6.7}) we  obtain (\ref{6.1}) .\rule{5pt}{5pt}
\bigskip
\begin{predl}\label{p6.2}
$q$-BMF $\makdo$ for $\nu>3/2$ can be represented as the
$q$-integral
$$
\makdo=\frac12 q^{-\nu^2+\nu}(1+q)
\Gamma_{q^2}(\nu+1)\sqrt{\frac{a_\nu}{a_{-\nu}}}\times
$$
\begin{equation}\label{6.8}
\times\int_0^\infty d_qz\frac{(-q^{2\nu+2}z^2,q^2)_\infty}
{(-z^2,q^2)_\infty}z
\phantom._0\Phi_3(-;0,0,q^2;q^2,-(\qpop q^{3/2}zs)^2)(s/2)^{-\nu}.
\end{equation}
\end{predl}

{\bf Proof.} Consider as before  the absolutely convergent $q$-integral
\begin{equation}\label{6.8b}
S_3^{(2)}(s)=\int_0^\infty d_qzg_\nu^{(2)}(z)
\phantom._0\Phi_3(-;0,0,q^2;q^2,-(\qpop q^{3/2}zs)^2).
\end{equation}
 It means according to Definition \ref{d2.1} and
(\ref{2.16}) that
$$
\lim_{m\to\pm\infty}q^m|g_\nu^{(2)}(q^m)\frac{2\p_s}{1+q}
J_0^{(2)}((1-q^2)q^{m-1}s,q^2)|=0
$$
It follows from (\ref{2.26}) that
\begin{equation}\label{6.8a}
\lim_{m\to\pm\infty}q^{2m}|g_\nu^{(2)}(q^m)
J_1^{(2)}((1-q^2)q^ms)|=0.
\end{equation}

Substitute $S_3^{(2)}(s)(s/2)^{-\nu}$ in
(\ref{4.2}). Then $S_3^{(2)}(s)$ satisfies  (\ref{5.10a}).
 Acting as
in in the proof of
Proposition \ref{p6.1} and using (\ref{2.16}), (\ref{2.25}), (\ref{2.26}),
  (\ref{2.27}) we come to the equality
$$
-q^{2\nu}\int_0^\infty d_qzg_\nu^{(2)}(z)z^{2\nu+1}\frac{2z^{-1}}
{1-q^2}[z^{-2\nu+1}\ddag J_1^{(2)}((1-q^2)zs;q^2)\ddag -
q^{-2\nu+1}z^{-2\nu+1}\ddag J_1^{(2)}((1-q^2)qzs;q^2)\ddag ]=
$$
$$
=\int_0^\infty d_qzg_\nu^{(2)}(z)(z/2)^{-1}\frac{2\p_z}
{1+q}(z/2\ddag J_1^{(2)}((1-q^2)zs;q^2)\ddag).
$$
Then by means of (\ref{2.3}) we obtain
$$
\lim_{m\to\infty}[g_\nu^{(2)}(q^{-m})(q^{2\nu-2m}+1)
J_1^{(2)}((1-q^2)q^{-m}s;q^2)-
g_\nu^{(2)}(q^m)(q^{2\nu+2m}+1)J_1^{(2)}((1-q^2)q^ms;q^2)]=
$$
\begin{equation}\label{6.9}
\int_0^\infty d_qz[\p_z(g_\nu^{(2)}(z)(z/2)^{-1})qz/2-
\p_z(g_\nu^{(2)}z^{2\nu+1})qz^{-2\nu+1}]
\ddag J_1^{(2)}((1-q^2)qzs;q^2))\ddag .
\end{equation}
It follows from (\ref{6.8a}) that the  left hand side of (\ref{6.9}) is equal to
zero. Then we have
\begin{equation}\label{6.10}
z^2[qg_\nu^{(2)}(z)-q^{2\nu+2}g_\nu^{(2)}(qz)]=
-qg_\nu^{(2)}(z)+g_\nu^{(2)}(qz)
\end{equation}
Note that
(\ref{6.10}) is the equation (\ref{3.4}) with $a=-q^{2\nu+2}, b=1,
\gamma=1$ and hence
\begin{equation}\label{6.11}
g_\nu^{(2)}(z)=\frac{(-q^{2\nu+2}z^2,q^2)_\infty}
{(-z^2,q^2)_\infty}z.
\end{equation}

Again $S_3^{(2)}(s)(s/2)^{-\nu}$ as a solution of (\ref{4.2}) can be
 represented in  the form
$$
\int_0^\infty d_qz\frac{(-q^{2\nu+2}z^2,q^2)_\infty}
{(-z^2,q^2)_\infty}z
\phantom._0\Phi_3(-;0,0,q^2;q^2,-(\qpop q^{3/2}zs)^2)(s/2)^{-\nu}=
$$
\begin{equation}\label{6.12}
=A\mbess+B\makdo.
\end{equation}
It can be derived from (\ref{5.a}) and (\ref{5.b}),
that $\lim_{s\to\infty}\mbess=\infty$ and \\
$\lim_{s\to\infty}\makdo=0$. Because the left hand side of (\ref{6.12})
vanishes if $s\to\infty,\quad A=0$. Multiplying (\ref{6.12}) on
$(s/2)^\nu$ and assuming $s=0$ we get
$$
\int_0^\infty d_qz\frac{(-q^{2\nu+2}z^2,q^2)_\infty}
{(-z^2,q^2)_\infty}z=B\frac{q^{-\nu^2+1/2}}
{4a_{-\nu}\sqrt{a_\nu a_{-\nu}}\sin{\nu\pi}\Gamma_{q^2}(-\nu+1)}.
$$
Using the same recipe as in  Proposition \ref{p6.1} we calculate
this $q$-integral:
$$
\int_0^\infty d_qz\frac{(-q^{2\nu+2}z^2,q^2)_\infty}
{(-z^2,q^2)_\infty}z=\frac{1-q}{1-q^{2\nu}}.
$$
Thus
$$
B=\frac{1-q}{1-q^{2\nu}}4q^{\nu^2-1/2}a_{-\nu}\sqrt{a_\nu a_{-\nu}}
\sin{\nu\pi}\Gamma_{q^2}(-\nu+1)=
\frac{2q^{\nu^2-\nu}}{(1+q)\Gamma_{q^2}(\nu+1)}
\sqrt{\frac{a_{-\nu}}{a_\nu}}
$$
and we obtain (\ref{6.8}) from (\ref{6.12}).\rule{5pt}{5pt}
\bigskip
\begin{rem}\label{r6.1}
If $q\to 1-0$ the equations (\ref{6.4}) and (\ref{6.10}) take the
form of the differential
equation
$$
z(1+z^2)g_\nu'(z)-[1-(2\nu+1)z^2]g_\nu(z)=0.
$$
The solution to this equation is
$$
g_\nu(z)=Cz(1+z^2)^{-\nu-1},
$$
and therefore we come to the classical integral representation of
Bessel-Macdonald function \cite{BE}
$$
K_\nu(s)=\Gamma(\nu+1)(s/2)^{-\nu}\int_0^\infty(1+z^2)^{-\nu-1}
zJ_0(zs)dz
$$
\end{rem}

\section{The representation of  $q$-BMF by the double $q$-integral}
\setcounter{equation}{0}
Here we combine  results from Sections 4 and 6.

Introduce new variable $\zeta$  with the commutations
\begin{equation}\label{7.1}
\zeta z=z\zeta,\qquad \zeta s=qs\zeta
\end{equation}
\begin{predl}\label{p7.1}
The $q$-BMF $\makd$ for $\nu>1/2$ can be represented as the double
$q$-integral
$$
\makd=\frac{q^{-\nu^2-\nu}(1+q)^2\Gamma_{q^2}(\nu+1)}
{4\phantom._2\Phi_1(q,q;q^3;q^2,1)}\sqrt{\frac{a_\nu}{a_{-\nu}}}
\times
$$
\begin{equation}\label{7.2}
\times\int_0^\infty d_qz\int_{-1}^1d_q\zeta
\frac{(-q^2z^2,q\zeta^2,q^2)_\infty}
{(-q^{-2\nu}z^2,\zeta^2,q^2)_\infty}z
\phantom._0\Phi_1(-;0;q,-i\qpop\zeta zs)(s/2)^{-\nu}
\end{equation}
\end{predl}

{\bf Proof.} This result follows from (\ref{6.1}) and (\ref{4.12}).
In view of of (\ref{7.1}) the result of $q$-integrations
is independent on their order
$$
\makd=\frac{q^{-\nu^2-\nu}(1+q)^2\Gamma_{q^2}(\nu+1)}
{4\phantom._2\Phi_1(q,q;q^3;q^2,1)}\sqrt{\frac{a_\nu}{a_{-\nu}}}
\times
$$
$$
\times\int_0^\infty d_qz\int_{-1}^1d_q\zeta
\frac{(-q^2z^2,q\zeta^2,q^2)_\infty}
{(-q^{-2\nu}z^2,\zeta^2,q^2)_\infty}z
E_q(-i\qpop(:(\zeta z)s:))(s/2)^{-\nu}=
$$
$$
=\frac{q^{-\nu^2-\nu}(1-q^2)^2\Gamma_{q^2}(\nu+1)}
{2\phantom._2\Phi_1(q,q;q^3;q^2,1)}\sqrt{\frac{a_\nu}{a_{-\nu}}}
\times
$$
$$
\times\sum_{m=-\infty}^\infty q^{2m}\sum_{l=0}^\infty
\frac{(-q^{2m+2},q^{2l+1},q^2)_\infty}
{(-q^{2m-2\nu},q^{2l},q^2)_\infty}\Cos_q(\qpop q^{m+l}s)(s/2)^{-\nu}.
$$
Inner series converges uniformly with respect to $m$, and we can change the order
of summation.\rule{5pt}{5pt}
\bigskip
\begin{rem}\label{r7.1}
If $q\to 1-0$  representation (\ref{7.2})
takes the form
\begin{equation}\label{7.3}
K_\nu(s)=\frac1\pi\Gamma(\nu+1)(s/2)^{-\nu}\int_0^\infty
(1+z^2)^{-\nu-1}zdz\int_{-1}^1(1-\zeta^2)^{-\frac12}
e^{-i\zeta zs}d\zeta.
\end{equation}
\end{rem}

Assume $\zeta=\cos\phi$ in (\ref{7.3}). Then
\begin{equation}\label{7.4}
K_\nu(s)=\frac1\pi\Gamma(\nu+1)(s/2)^{-\nu}\int_0^\infty
(1+z^2)^{-\nu-1}zdz\int_0^\pi e^{-izs\cos\phi}d\phi
\end{equation}
It is known that  BMF is the two-dimensional Fourier transform of
 $f(x,y)=(1+x^2+y^2)^{-\nu-1}$:
$$
\int\int(1+x^2+y^2)^{-\nu-1}e^{-i(x\xi+y\eta)}dxdy=
\frac{2^{-\nu}\pi}{\Gamma(\nu+1)}K_\nu(\sqrt{\xi^2+\eta^2})
(\xi^2+\eta^2)^{\nu/2}.
$$
The change of variables $x=z\cos\alpha, y=z\sin\alpha,\quad
\xi=s\cos\beta, \eta=s\sin\beta$, and $\alpha-\beta=\phi$
leads  to (\ref{7.4}).

Thus  (\ref{7.2})  can be considered as the $q$-analog of
the Fourier transform  of
$$
f(z,\zeta)=\frac{(-q^2z^2,q\zeta^2,q^2)_\infty}
{(-q^{-2\nu}z^2,\zeta^2,q^2)_\infty}
$$
on the quantum plane  in the polar coordinate
representation (\ref{7.1}).

\section{The special case of commuting variables}
\setcounter{equation}{0}

In this section we  assume that
\begin{equation}\label{8.1}
zs=sz,\qquad \zeta s=s\zeta,\qquad z\zeta=\zeta z.
\end{equation}
Then we have the $q$-integral representations which we will give
without proof.
\begin{equation}\label{8.2}
\mbes=\frac{(1+q)(s/2)^\nu}
{2\Gamma_{q^2}(\nu+1/2)\Gamma_{q^2}(1/2)}
\int_{-1}^1\frac{(q^2z^2,q^2)_\infty}
{(q^{2\nu+1}z^2,q^2)_\infty}e_q(\qpop zs)d_qz,
\end{equation}
$$
\mbess=
$$
\begin{equation}\label{8.3}
=\frac{(s/2)^\nu}
{2\Gamma_{q^2}(\nu+1)\phantom._2\Phi_1(q^{-2\nu+1},q;q^3;q^2,1)}
\int_{-1}^1\frac{(q^{-2\nu+1}z^2,q^2)_\infty}
{(z^2,q^2)_\infty}E_q(\qpop zs)d_qz,
\end{equation}
$$
\makd=
$$
\begin{equation}\label{8.4}
=\frac{q^{-\nu^2+1/2}\Gamma_{q^2}(\nu+1/2)\Gamma_{q^2}(1/2)}
{4{\bf Q_\nu}}\sqrt{\frac{a_\nu}{a_{-\nu}}}(s/2)^{-\nu}
\int_{-\infty}^\infty\frac{(-q^2z^2,q^2)_\infty}
{(-q^{-2\nu+1}z^2,q^2)_\infty}e_q(i\qpop zs)d_qz,
\end{equation}
$$
\makdo=
$$
\begin{equation}\label{8.5}
=\frac{q^{-\nu^2+\nu}\Gamma_{q^2}(\nu+1/2)\Gamma_{q^2}(1/2)}
{4{\bf Q_{1/2}}}\sqrt{\frac{a_\nu}{a_{-\nu}}}(s/2)^{-\nu}
\int_{-\infty}^\infty\frac{(-q^{2\nu+1}z^2,q^2)_\infty}
{(-z^2,q^2)_\infty}E_q(i\qpop zs)d_qz,
\end{equation}
$$
\makd=
$$
\begin{equation}\label{8.6}
=\frac12q^{-\nu^2-\nu}(1+q)\Gamma_{q^2}(\nu+1)
\sqrt{\frac{a_\nu}{a_{-\nu}}}(s/2)^{-\nu}
\int_0^\infty\frac{(-q^2z^2,q^2)_\infty}
{(-q^{-2\nu+1}z^2,q^2)_\infty}zJ_0^{(1)}((1-q^2)zs;q^2)d_qz,
\end{equation}
$$
\makdo=
$$
\begin{equation}\label{8.7}
=\frac12q^{-\nu^2+\nu}(1+q)\Gamma_{q^2}(\nu+1)
\sqrt{\frac{a_\nu}{a_{-\nu}}}(s/2)^{-\nu}
\int_0^\infty\frac{(-q^{2\nu+2}z^2,q^2)_\infty}
{(-z^2,q^2)_\infty}zJ_0^{(2)}((1-q^2)zs;q^2)d_qz,
\end{equation}
$$
\makd=\frac{q^{-\nu^2-\nu}(1+q)^2\Gamma_{q^2}(\nu+1)}
{4\Gamma_{q^2}^2(1/2)}\sqrt{\frac{a_\nu}{a_{-\nu}}}(s/2)^{-\nu}
\times
$$
\begin{equation}\label{8.8}
\times\int_0^\infty zd_qz\int_{-1}^1
\frac{(-q^2z^2,q^2\zeta^2,q^2)_\infty}
{(-q^{-2\nu}z^2,q\zeta^2,q^2)_\infty}e_q(-i\qpop z\zeta s)d_q\zeta,
\end{equation}
$$
\makdo=\frac{q^{-\nu^2+\nu}(1+q)^2\Gamma_{q^2}(\nu+1)}
{4\phantom._2\Phi_1(q,q;q^3;q^2,1)}
\sqrt{\frac{a_\nu}{a_{-\nu}}}(s/2)^{-\nu}
\times
$$
\begin{equation}\label{8.9}
\times\int_0^\infty zd_qz\int_{-1}^1
\frac{(-q^{2\nu+2}z^2,q\zeta^2,q^2)_\infty}
{(-z^2,\zeta^2,q^2)_\infty}E_q(-i\qpop z\zeta s)d_q\zeta,
\end{equation}

{\bf Acknowledgments.}
{\sl M.O is
 grateful to  the Max-Planck Institute for Mathematics in Bonn,
where this work was completed.}

\small{

}


\begin{thebibliography}{40}
\bibitem{OR1}\  Olshanetsky M. and Rogov V., {\em The modified
$q$-Bessel and the $q$-Bessel-Macdonald Functions},
 Preprint ITEP-TH-6/95, q-alg/950913
\bibitem{Ja}\ Jackson F.H. {\em The application of basic
numbers to Bessel's and
Legendre's functions}, Proc. London math. Soc. (2) {\bf 2} (1905) 192-220
\bibitem{Vi}\ Vilenkin, N.Ya., Special Functions and the Theory of Group
Representations, Amer. Math. Soc. Transl. of Math. Monographs {\bf 22},
Providence 1968
\bibitem{VK1}\ Vaksman L. and Korogodskii L., {\em Algebra of bounded functions
on the quantum group of the motions of the plane and q-analogues of Bessel
functions}, Soviet. Math. Dokl. {\bf 39} (1989) 173-177
\bibitem{FV}\ Floreanini R. and Vinet L. {\em Addition formulas
for q-Bessel functions},
preprint University of Montreal, UdeM-LPN-TH60, (1991);
{\em Representations of quantum algebras
and q-special functions}, preprint Universite de Montreal (1991)
\bibitem{Kol}\ Koelink H.T., {\em The quantum group of plane motions
and the Hahn-Exton $q$-Bessel function},  Duke Math. J., {\bf 76}
(1994) 483-508
\bibitem{J}\ Jacquet H. P. B., ~{\em Fonctions de Whittaker associees aux
 groupes de
Chevalley},~Bull.Soc. Math. France, 95 (1967) 243-309
\bibitem{Sh}\ Schiffmann G.,
 {\em Integrales d'entrelacement
 et fonctions de Whittaker},~Bull.Soc.Math.France, 99 (1971) 3-72
\bibitem{Ko} B.Kostant, {\em Whittaker vectors and
representation theory},  Invent. Math., {\bf 48} (1978) 101-184
\bibitem{Ha} M.Hashizume,~{\em Whittaker functions on semisimple
 Lie groups},
 Hiroshima Math. Journ., {\bf 12} (1982) 259-293
\bibitem{STS}\ Semenov-Tian-Shansky M.,
~~ {\em Quantization of the open Toda chains}, ~in
"Sovremenie problemi matematiki" ~VINITI, v.16 (1987) 194-226, (in Russian)
\bibitem{OR2}\ Olshanetsky M. and Rogov V.,
{\em Liouville quantum mechanics on a
lattice from geometry of quantum Lorentz group }, Journ. Phys. A: Math. Gen,
{\bf 27} (1994) 4669-4683
\bibitem{CZ}\ Chryssomalakos Ch. and Zumino B., {\em Translations,
Integrals and
Fourier Transforms in the Quantum Plane}, Salamfest, 1993, 327-346
\bibitem{VK}\ Vaksman L. and Korogodsky L.,
 {\em Harmonic Analisys on Quantum
Hyperboloids}, Preprint ITP-90-27P (in Russian)
 \bibitem{GR}\ Gasper G. and Rahman M.,
~{\em Basic Hypergeometric Series}, Cambridge:
 Cambridge University Press, (1990)
\bibitem{BMP} A.Brychkov, Yu.Prudnikov and O.Marychev,
{\em Integrals and Series,} vol.1, Nauka, Moscow 1986
\bibitem{BE}\ Bateman H., and Erdlyi A., {\em Higher transcendental functions},
 v.2
Mc Graw-Hill Book Company, 1953

\end{thebibliography}
\end{document}